\documentclass[12pt]{amsart}
\usepackage{amssymb}

\expandafter\ifx\csname delta.sty\endcsname\relax \else\endinput\fi
\expandafter\edef\csname delta.sty\endcsname{%
 \catcode`\noexpand\@=\the\catcode`\@\space}

\let\atbefore @

\catcode`\@=11

\newif\ifMag
\let\@ft@\expandafter
\numberwithin{equation}{section}

\newif\ifComments

\def\tod@y{\ifcase\month\or
 January\or February\or March\or April\or May\or June\or July\or
 August\or September\or October\or November\or December\fi\space\,
\number\day,\space\,\number\year}

\def\h@@r{hh}\def\m@n@te{mm}
\def\wh@tt@me{\count@\time\divide\count@ 60\edef\h@@r{\number\count@}%
 \multiply\count@ -60\advance\count@\time\edef
 \m@n@te{\ifnum\count@<10 0\fi\number\count@}}
\def\t@me{\h@@r\/{\rm:}\m@n@te} \let\whattime\wh@tt@me
\let\Today\tod@y \let\nowtime\t@me

\def\ftext#1{{\let\thefootnote\relax\footnotetext{\vsk-.8>\nt #1}}}
\def\em#1{{\itshape #1\/}}

 {\end{trivlist}\egroup\global\@ignoretrue\ifComments\unvbox\commentbox\fi
 \noindent}

\def\gadv{\global\adv} \def\gad#1{\gadv#1\@ne} \def\gadneg#1{\gadv#1-\@ne}
\def\textindent#1{\indent\llap{#1\enspace}\ignorespaces}

\newcount\itemlet
\def\newbi{\itemlet 96} \newbi
\def\bitem{\gad\itemlet\endgraf\hangindent1.5\parindent
 \hglue-.5\parindent\textindent{\upshape\rlap{\char\the\itemlet}\hp{b})}}

\newcount\itemrm

\def\iitem{\gad\itemrm\endgraf\hangindent1.5\parindent\hglue-.5\parindent
 \textindent{\upshape\hp{v}\llap{\romannumeral\the\itemrm})}}

\let\Disp\[ \let\endD\] \let\{\protect

\def\Tag#1{\label{e:#1}\let\notag\relax} 
\def\sh@nd#1#2{\begin{#1*}#2\end{#1*}}
\def\n@t@gs#1#2#3{\let\n@@@l\\ \begin{#1}#2\global\let\d@bl\\
 \gdef\\{\notag\d@bl}#3\notag\global\let\\\d@bl\end{#1}\let\\\n@@@l}

\def\Gather#1\endG{\n@t@gs{gather}{}{#1}}
\def\gAther#1\endG{\sh@nd{gather}{#1}}
\def\Align#1\endA{\n@t@gs{align}{}{#1}}
\def\aLign#1\endA{\sh@nd{align}{#1}}
\def\Alignat#1#2\endAt{\n@t@gs{alignat}{#1}{#2}}
\def\aLignat#1\endAt{\sh@nd{alignat}{#1}}

\def\(#1){\textup{(\ref{e:#1})}}
\def\[{\@ifnextchar:\c@t@sect\c@t@}
\def\c@t@sect:#1]{\ref{s:#1}} \def\c@t@#1]{\ref{t:#1}}

\def\qed{\hbox{}\nobreak\hfill\nobreak{\m@th$\,\square$}}
\def\sk@@p#1{\par\skip@#1\relax\ifdim\lastskip<\skip@\relax\removelastskip
 \vskip\skip@\fi}
\def\demo#1{\sk@@p\medskipamount\nt{\ignore\it #1\unskip.}\enspace
 \ignore}
\def\enddemo{\sk@@p\medskipamount}
 
\def\Pf#1.{\demo{Proof #1}}

\let\bls\baselineskip \let\ignore\ignorespaces \let\adv\advance
\def\vsk#1>{\vskip#1\bls}
\def\vv#1>{\vadjust{\vsk#1>}\ignore}
\def\vvn#1>{\vadjust{\nobreak\vsk#1>\nobreak}\ignore}
\def\vvv#1>{\vskip\z@\vsk#1>\nt\ignore}
\def\vvgood{\vadjust{\penalty-500}} 

\def\mathbox#1{\hbox{\m@th$#1$}}

\let\dsize\displaystyle \let\tsize\textstyle
\let\ssize\scriptstyle \let\sss\scriptscriptstyle
 
\let\vp\vphantom \let\hp\hphantom \let\nt\noindent
\def\hline{\hbox to\hsize}
\let\cline\centerline \let\lline\leftline \let\rline\rightline
\def\nn#1>{\noalign{\vskip#1\p@}} \def\NN#1>{\openup#1\p@}
 
\let\Lim\lim \def\lim{\Lim\limits} \let\Sum\sum \def\sum{\Sum\limits}
 
\let\Prod\prod \def\prod{\Prod\limits} \let\Int\int \def\int{\Int\limits}

\def\tsum{\mathop{\tsize\Sum}\limits} 

\def\~{\leavevmode\@ifnextchar~\m@n@s\@md@sh}
\def\m@n@s~{\raise.15ex\mathbox{-}} \def\@md@sh{\raise.13ex\hbox{--}}
\let\procent\% \def\%#1{\ifmmode\mathop{#1}\limits\else\procent#1\fi}
\let\@ml@t\" \def\"#1{\ifmmode ^{(#1)}\else\@ml@t#1\fi}
\let\@c@t@\' \def\'#1{\ifmmode _{(#1)}\else\@c@t@#1\fi}
\let\colon\: \def\:{^{\vp|}} \def\&{.\kern.1em} \def\^#1{\text{\m@th#1}}

\newif\ifNewskips

\def\Newskips{\global\Newskipstrue
 \gdef\>{\relax\ifmmode\mskip.666667\thinmuskip\relax\else\kern.111111em\fi}
 \gdef\}{\relax\ifmmode\mskip-.666667\thinmuskip\relax\else\kern-.111111em\fi}
 \gdef\){\relax\ifmmode\mskip.333333\thinmuskip\relax\else\kern.0555556em\fi}
 \gdef\]{\relax\ifmmode\mskip-.333333\thinmuskip\relax\else\kern-.0555556em\fi}}
\Newskips

\def\End{\mathop{\mathrm{End}\>}}

\def\id{\mathrm{id}}  
\def\1{^{-1}} \def\_#1{_{\Rlap{#1}}}
\def\vst#1{{\lower1.9\p@
 \hbox{\m@th$\bigr|_{\raise.5\p@\hbox{\m@th$\ssize#1$}}$}}}
\def\vrp#1:#2>{{\vrule height#1 depth#2 width\z@}}
\def\vru#1>{\vrp#1:\z@>} \def\vrd#1>{\vrp\z@:#1>}
 
\def\sscr#1{\raise.3ex\hbox{\m@th$\sss#1$}} \def\@@PS{\mathbf{OOPS!!!}}

\def\lsym#1{#1\alb\ldots\relax#1\alb}
\def\lc{\lsym,}   \def\lox{\lsym\ox}

\let\texspace\ \def\ {\ifmmode\alb\fi\texspace}

\def\Line#1{\kern-.5\hsize\hline{\m@th$\dsize#1$}\kern-.5\hsize}
\def\Lline#1{\kern-.5\hsize\lline{\m@th$\dsize#1$}\kern-.5\hsize}
\def\Cline#1{\kern-.5\hsize\cline{\m@th$\dsize#1$}\kern-.5\hsize}
\def\Rline#1{\kern-.5\hsize\rline{\m@th$\dsize#1$}\kern-.5\hsize}

\def\Ll@p#1{\llap{\m@th$#1$}} \def\Rl@p#1{\rlap{\m@th$#1$}}
 \def\Cl@p#1{\llap{\m@th$#1$\hss}}
\def\Llap#1{\mathchoice{\Ll@p{\dsize#1}}{\Ll@p{\tsize#1}}{\Ll@p{\ssize#1}}%
 {\Ll@p{\sss#1}}}
\def\Clap#1{\mathchoice{\Cl@p{\dsize#1}}{\Cl@p{\tsize#1}}{\Cl@p{\ssize#1}}%
 {\Cl@p{\sss#1}}}
\def\Rlap#1{\mathchoice{\Rl@p{\dsize#1}}{\Rl@p{\tsize#1}}{\Rl@p{\ssize#1}}%
 {\Rl@p{\sss#1}}}
 
\def\LRtph#1#2{\setbox\z@\hbox{#1}\dimen\z@\wd\z@\hbox{\hbox to\dimen\z@{#2}}}
\def\LRph#1#2{\LRtph{\m@th$#1$}{\m@th$#2$}}

\def\Lto#1{\setbox\z@\hbox{\m@th$\tsize{#1}$}%
 \mathrel{\mathop{\hbox to\wd\z@{\rightarrowfill}}\limits#1}}
\def\Lgets#1{\setbox\z@\hbox{\m@th$\tsize{#1}$}%
 \mathrel{\mathop{\hbox to\wd\z@{\leftarrowfill}}\limits#1}}

\def\vpb#1{{\vp{\big(}}^{\]#1}}

\let\alb\allowbreak

 \let\x\times \let\ox\otimes 
 
\let\le\leqslant 
\let\der\partial \let\8\infty \let\*\star

\let\map\mapsto

\let\lb\lbrace \let\rb\rbrace

\let\Bbb\mathbb

\let\frak\mathfrak

\def\pms{\raise.25ex\mathbox{\ssize\pm}\>}
\def\mps{\raise.25ex\mathbox{\ssize\mp}\>}

\let\dl\delta \let\Dl\Delta 
\let\epe\epsilon \let\eps\varepsilon \let\epsilon\eps

\let\la\lambda

 \let\phi\varphi

\def\C{\Bbb C}
\def\R{\Bbb R}

\def\h@ph{\discretionary{}{}{-}} \def\$#1$-{\,\text{\m@th$#1$}\h@ph}

\def\difl/{differential} \def\dif/{difference}
\def\cf.{cf.\ \ignore} \def\Cf.{Cf.\ \ignore}
\def\egv/{eigenvector} \def\eva/{eigenvalue} \def\eq/{equation}
\def\lhs/{the left hand side} \def\rhs/{the right hand side}
\def\Lhs/{The left hand side} \def\Rhs/{The right hand side}
\def\gby/{generated by} \def\wrt/{with respect to} \def\st/{such that}
\def\resp/{respectively} \def\off/{offdiagonal} \def\wt/{weight}
\def\pol/{polynomial} \def\rat/{rational} \def\tri/{trigonometric}
\def\fn/{function} \def\var/{variable} \def\raf/{\rat/ \fn/}
\def\inv/{invariant} \def\hol/{holomorphic} \def\hof/{\hol/ \fn/}
\def\mer/{meromorphic} \def\mef/{\mer/ \fn/} \def\mult/{multiplicity}
\def\sym/{symmetric} \def\perm/{permutation}
\def\rep/{representation} \def\irr/{irreducible} \def\irrep/{\irr/ \rep/}
\def\hom/{homomorphism} \def\aut/{automorphism} \def\iso/{isomorphism}
\def\lex/{lexicographical} \def\as/{asymptotic} \def\asex/{\as/ expansion}
\def\ndeg/{nondegenerate} \def\neib/{neighbourhood} \def\deq/{\dif/ \eq/}
\def\hw/{highest \wt/} \def\gv/{generating vector} \def\eqv/{equivalent}
\def\msd/{method of steepest descend} \def\pd/{pairwise distinct}
\def\wlg/{without loss of generality} \def\Wlg/{Without loss of generality}
\def\onedim/{one-dim\-en\-sion\-al} \def\fd/{fi\-ni\-te-dim\-en\-sion\-al}
\def\qcl/{quasiclassical} \def\hwv/{\hw/ vector}
\def\hgeom/{hyper\-geo\-met\-ric} \def\hint/{\hgeom/ integral}
\def\hwm/{\hw/ module} \def\emod/{evaluation module} \def\Vmod/{Verma module}
\def\symg/{\sym/ group} \def\sol/{solution} \def\eval/{evaluation}
\def\anf/{analytic \fn/} \def\anco/{analytic continuation}
\def\qg/{quantum group} \def\qaff/{quantum affine algebra}

\def\Rm/{\^{$R$-}matrix} \def\Rms/{\^{$R$-}matrices}
\def\YB/{Yang-Baxter \eq/}
\def\Ba/{Bethe ansatz} \def\Bv/{Bethe vector} \def\Bae/{\Ba/ \eq/}
\def\KZv/{Knizh\-nik-Zamo\-lod\-chi\-kov} \def\KZvB/{\KZv/-Bernard}
\def\KZ/{{\sl KZ\/}} \def\qKZ/{{\sl qKZ\/}}
\def\KZB/{{\sl KZB\/}} \def\qKZB/{{\sl qKZB\/}}
\def\qKZo/{\qKZ/ operator} \def\qKZc/{\qKZ/ connection}
\def\KZe/{\KZ/ \eq/} \def\qKZe/{\qKZ/ \eq/} \def\qKZBe/{\qKZB/ \eq/}

\def\LPT/{Laboratoire de Physique Th\'eorique ENSLAPP}
\def\ENSLyon/{\'Ecole Normale Sup\'erieure de Lyon}
\def\DMS/{Department of Mathematics, Faculty of Science}
\def\DMO/{\DMS/, Osaka University}
\def\DMOaddr/{Toyonaka, Osaka 560, Japan}
\def\dmoemail/{vt@math.sci.osaka-u.ac.jp}
\def\SPb/{St\&Petersburg}
\def\home/{\SPb/ Branch of Steklov Mathematical Institute}
\def\homeaddr/{Fontanka 27, \SPb/ \,191011, Russia}
\def\homemail/{vt@pdmi.ras.ru}
\def\absence/{On leave of absence from \home/}
\def\UNC/{Department of Mathematics, University of North Carolina}
\def\ChH/{Chapel Hill}
\def\UNCaddr/{\ChH/, NC 27599, USA} \def\avemail/{av@math.unc.edu}
\def\grant/{NSF grant DMS\~9501290}	
\def\Grant/{Supported in part by \grant/}

\def\Aomoto/{K\&Aomoto}
\def\Dri/{V\]\&G\&Drin\-feld}
\def\Fadd/{L\&D\&Fad\-deev}
\def\Feld/{G\&Felder}
\def\Fre/{I\&B\&Fren\-kel}
\def\Gustaf/{R\&A\&Gustafson}
\def\Kazh/{D\&Kazhdan} \def\Kir/{A\&N\&Kiril\-lov}
\def\Kor/{V\]\&E\&Kore\-pin}
\def\Lusz/{G\&Lusztig}
\def\MN/{M\&Naza\-rov}
\def\Resh/{N\&Reshe\-ti\-khin} \def\Reshy/{N\&\]Yu\&Reshe\-ti\-khin}
\def\Skl/{E\&K\&Sklya\-nin}
\def\SchV/{V\]\&\]V\]\&Schecht\-man} \def\Sch/{V\]\&Schecht\-man}
\def\Takh/{L\&A\&Takh\-tajan}
\def\VT/{V\]\&Ta\-ra\-sov} \def\VoT/{V\]\&O\&Ta\-ra\-sov}
\def\Varch/{A\&\]Var\-chenko} \def\Varn/{A\&N\&\]Var\-chenko}

\def\AMS/{Amer.\ Math.\ Society}
\def\CMP/{Comm.\ Math.\ Phys.{}}
\def\DMJ/{Duke.\ Math.\ J.{}}
\def\Inv/{Invent.\ Math.{}} 
\def\IMRN/{Int.\ Math.\ Res.\ Notices}
\def\JPA/{J.\ Phys.\ A{}}
\def\JSM/{J.\ Soviet\ Math.{}}
\def\LMP/{Lett.\ Math.\ Phys.{}}
\def\LMJ/{Leningrad Math.\ J.{}}
\def\LpMJ/{\SPb/ Math.\ J.{}}
\def\SIAM/{SIAM J.\ Math.\ Anal.{}}
\def\SMNS/{Selecta Math., New Series}
\def\TMP/{Theor.\ Math.\ Phys.{}}
\def\ZNS/{Zap.\ nauch.\ semin. LOMI}

\def\ASMP/{Advanced Series in Math.\ Phys.{}}

\def\AMSa/{AMS \publaddr Providence}
\def\Birk/{Birkh\"auser}
\def\CUP/{Cambridge University Press} \def\CUPa/{\CUP/ \publaddr Cambridge}
\def\Spri/{Springer-Verlag} \def\Spria/{\Spri/ \publaddr Berlin}
\def\WS/{World Scientific} \def\WSa/{\WS/ \publaddr Singapore}

\csname delta.sty\endcsname

\newif\ifUS
\ifx\USversion\relax\UStrue\fi

\ifUS
\fi

\evensidemargin\oddsidemargin
\edef\restoreatcode{\catcode`\noexpand\@=\the\catcode`\@}
\makeatletter

\newtheorem{thm}{Theorem}

\newenvironment{abstr}{\begingroup\narrower\small
\nt{\sc Abstract.}\enspace\ignorespaces}{\endgraf\endgroup}

\def\fratop{\genfrac{}{}{\z@}1}

\def\Ref#1{{\rm(\ref{#1})}}
\def\no.{no.\,\,\ignore}

\def\dash{\raise.13ex\hbox{-}}

\def\cnn#1>{\\[-\bls]\noalign{\vsk#1>}\notag}
\def\mmgood#1:#2>{\\\noalign{\vsk#1>\penalty-500\vsk#2>}\notag}

\def\fratop{\genfrac{}{}{0pt}1}
\def\q-{\>\hbox{\m@th$q\)$-}}

\def\gl{\frak{gl}} \def\gsl{\frak{sl}}
\def\gln{\gl_N}

\def\zn{z_1\lc z_n} \def\lan{\la_1\lc\la_N}

\def\Vox{V_1\lox V_n}

\def\R-{\>\hbox{\m@th$R\)$-}}

\newcommand{\bean}{\begin{eqnarray}}
\newcommand{\eean}{\end{eqnarray}}
\newcommand{\bea}{\begin{eqnarray*}}
\newcommand{\eea}{\end{eqnarray*}}

\restoreatcode

\begin{document}

\title[Dynamical Differential Equations]
{Dynamical Differential Equations\\[3pt]
Compatible with Rational qKZ Equations}

\author[\smash{V\]\&Tarasov and A\&\]Varchenko}]
{\vbox{}V\]\&Tarasov$^\star$ \>and \;A\&\]Varchenko$^\diamond$}

\ftext{\mathsurround 0pt
$\]^\star\)$Supported in part by RFFI grant 02\)\~\)01\~\)00085a
\>and \,CRDF grant RM1\~\)2334\)\~MO\)\~\)02\vv.2>\\
\hp{$^*$}{\normalsize\sl E-mail\/{\rm:} \homemail/\,{\rm,}
vtarasov@math.iupui.edu}\vv.1>\\
${\]^\diamond\)}$Supported in part by NSF grant DMS\)\~\)0244579\vv.2>\\
\vv-1.2>
\hp{$^*$}{\normalsize\sl E-mail\/{\rm:} anv@email.unc.edu}}

\hrule width0pt

\maketitle

\begin{center}
{\it\=$^\star$\home/\\[2pt]\homeaddr/\\[6pt]
$^\star$Department of Mathematical Sciences,\\[2pt]
Indiana University Purdue University at Indianapolis,\\[2pt]
Indianapolis, IN 46202, USA\\[6pt]
$^\diamond$\UNC/ at \ChH/\\[2pt]\UNCaddr/}

\vsk1.5>
{\sl March 2004}
\end{center}

\vsk1.5>
\begin{abstr}
For the Lie algebra $\gln\}$ we introduce a system of differential operators
called the dynamical operators. We prove that the dynamical differential
operators commute with the $\gln\}$ rational quantized Knizhnik-Zamolodchikov
difference operators. We describe the transformations of the dynamical
operators under the natural action of the $\gln\}$ Weyl group.
\end{abstr}

\thispagestyle{empty}

\vsk>
\vsk0>
\section{Introduction}

The rational quantized Knizhnik-Zamolodchikov difference equations associated
with the Lie algebra $\gln\}$ is a system of difference equations of the form
\vvn.3>
\begin{align}
\label{KZ}
U(z_1\lc z_i\]+p\lc z_n;{}& \lan)\,={}
\\[4pt]
{}=\,K_i(\zn;{}&\lan)\>U(\zn;\lan)\,,
\notag
\end{align}
$i=1\lc n$. Here $p$ is the step of the difference equations, $U(z\);\la)$
is a function with values in the tensor product $\Vox$ of $n$ highest weight
\$\gln$-modules, $K_i(z\);\la)$ is a suitable linear operator on the tensor
product. The \qKZ/ difference equations have found many applications; for
example, see~\cite{EFK}\), \cite{JM}\), \cite{V2}. In \cite{TV5} we suggested
a system of differential equations of the form
\vvn.3>
\begin{equation}
\label{DYN}
\Bigl(p\>\la_a\frac\der{\der\la_a}\,+\,L_a(\zn;\lan)\Bigr)\,
U(\zn;\lan)\,=\,0\,,
\vv.2>
\end{equation}
$a=1\lc N$. Here $L_a(z\);\la)$ is a suitable linear operator on ${\Vox}$.
\vvgood
We called this system the dynamical differential equations. In this paper we
prove that the \qKZ/ difference equations and dynamical differential equations
are compatible, and describe the transformation properties of the dynamical
equations under the natural action of the $\gln\}$ Weyl group.
\vsk.2>
The phenomenon of existence of dynamical equations compatible with \KZ/
equations was discovered in \cite{FMTV}. There are many versions of \KZ/
equations: rational and trigonometric, differential and difference. For each
version of the \KZ/ equations there exists a complementary system of dynamical
equations which is compatible with the \KZ/ equations. In \cite{FMTV}
the rational differential \KZ/ equations were considered and the compatible
dynamical differential equations were introduced. In \cite{TV4} the
trigonometric differential \KZ/ equations were considered and the compatible
dynamical difference equations were introduced. In \cite{EV} the trigonometric
difference \KZ/ equations were considered and the compatible dynamical
difference equations were introduced.
\vsk.1>
All versions of \KZ/ equations have hypergeometric solutions, see \cite{SV}\),
\cite{V1}\), \cite{TV1}\), \cite{TV2}\), \cite{TV3}\), \cite{FV}\),
\cite{FTV}\). The general conjecture is that the hypergeometric solutions also
satisfy the corresponding dynamical equations. For rational \KZ/ differential
equations that was proved in \cite{FMTV}\), for trigonometric \KZ/ differential
equations that was proved in \cite{MV}\). We plan to prove that the
hypergeometric solutions of the rational \qKZ/ difference equations \Ref{KZ}
also satisfy the dynamical differential equations \Ref{DYN} in our next paper.
\vsk.1>
The fact that the hypergeometric solutions of the \KZ/ equations satisfy also
the additional dynamical equations is useful for aplications. For example,
the dynamical equations, in principal, allow us to recover the hypergeometric
solution from the asymptotics of the solution as $\la$ tends to a special
value. In \cite {TV6} we used the dynamical equations in that way to find
a formula for Selberg type integrals associated with $\gsl_3$.
\vsk.1>
There is also another phenomenon: the \KZ/ and dynamical equations correspond
to each other under the $(\gln\>,\gl_n)$ duality. This phenomenon was
discovered in \cite{T} and \cite{TV5}\). It turns out that under the
$(\gln\>,\gl_n)$ duality the \KZ/ and dynamical equations associated with
$\gln$ become respectively the dynamical and \KZ/ equations associated with
$\gl_n$. That kind of the duality and the (partially conjectural) fact that
hypergeometric solutions satisfy both the \KZ/ and dynamical equations, in
principal, allows us to identify hypergeometric solutions of the $\gln\}$
and $\gl_n\}$ \KZ/ equations. We used that idea in \cite{TV7} and \cite{TV8}
to prove certain nontrivial identities between hypergeometric integrals of
different dimensions.

\vsk-.5>\vsk0>
\section{\qKZ/ and dynamical operators}

\vsk.3>
\subsection{The Yangian $Y(\gln)$ and the rational \R-matrix}

Let $e_{a,b}$, $a,b=1,\ldots N$, be the standard generators of the Lie algebra
$\gln\}$, $[ e_{a,b}\,, \, e_{c,d}]\,=\,\dl_{bc}\,e_{a,d}\,-\,\dl_{ad}
\,e_{c,b}\,. $
\vsk.1>
The Yangian $Y(\gln)$ is the unital associative algebra with
generators $\{T_{a,b}^{(s)}\}$ where $a\),b=1\lc N$ and $s=1,2,\ldots$.
Organize them into generating series
$$
T_{a,b}(u)\,=\,\dl_{a,b}\,+\sum_{s=1}^\infty T_{a,b}^{(s)}u^{-s}\,.
\vv-.1>
$$
The defining relations in $Y(\gln)$ have the form
\vvn.2>
\begin{equation}
\label{yan-rel}
\bigl[\)T_{a,b}(u)\>,T_{c,d}(v)\)\bigr]\,=\;
\frac{T_{c,b}(v)\>T_{a,d}(u)-T_{c,b}(u)\>T_{a,d}(v)}{u-v}\;.
\end{equation}
The Yangian $Y(\gln)$ is a Hopf algebra with coproduct
$\Dl:Y(\gln)\to Y(\gln)\ox Y(\gln)$,
\bea
\Dl\>:\>T_{a,b}(u)\,\map\,\sum_{c=1}^N\,T_{c,b}(u)\ox T_{a,c}(u)\,.
\eea
There is a one-parametric family of automorphism $\rho_x\}:Y(\gln)\to Y(\gln)$,
$$
\rho_x:\>T_{a,b}(u)\,\map\,T_{a,b}(u-x)\,.
$$
\par
The Yangian $Y(\gln)$ contains the universal enveloping algebra $U(\gln)$ as
a Hopf subalgebra. The embedding is defined by $e_{a,b}\map T_{b,a}^{(1)}$ for
all $a\),b=1\lc N$. We identify $U(\gln)$ with its image in $Y(\gln)$ under
this embedding.
\vsk.1>
There is an evaluation homomorphism $\epe:Y(\gln)\to U(\gln)$,
\begin{equation}
\label{evalu}
\epe\>:\>T_{a,b}(u)\,\map\,\dl_{a,b}\)+\)e_{b,a}\)u^{-1}\).
\end{equation}
Both the automorphism $\rho_x$ and the homomorphism $\epe$ restricted to
the subalgebra $U(\gln)$ are the identity maps.
\vsk.1>
For a \$\gln$-module $V$ denote by $V(x)$ the \$Y(\gln)$-module
induced from $V$ by the homomorphism $\epe\circ\rho_x$.
The module $V(x)$ is called an evaluation module.
\vsk.2>
Let $V_1, V_2$ be Verma modules over $\gln\}$ with highest weight vectors
$v_1, v_2$, respectively. For generic complex numbers $x, y$ the Yangian
modules $V_1(x)\ox V_2(y)$ and $V_2(y)\ox V_1(x)$ are known to be isomorphic.
An isomorphism of the modules sends $\C\)(v_1\ox v_2)$ to $\C\)(v_2\ox v_1)$.
We fix an isomorphism requiring that the vector $v_1\ox v_2$ is mapped to
$v_2\ox v_1$. The isomorphism has the form
$$
P\)R_{V_1,V_2}(x-y):V_1(x)\ox V_2(y)\,\to\,V_2(y)\ox V_1(x)
\vv.4>
$$
where ${P:V_1\ox V_2\to V_2\ox V_1}$ is the permutation of factors, and
$R_{V_1,V_2}(x)$ takes values in $\End(V_1\ox V_2)$. The operator
$R_{V_1,V_2}(x)$ respects the weight decomposition of the \$\gln$-module
$V_1\ox V_2$ and its restriction to any weight subspace is a rational function
of $x$. The operator $R_{V_1,V_2}(x)$ is called the rational \R-matrix for
the tensor product $V_1\ox V_2$.
\vsk.2>
The definition of $R_{V_1,V_2}(x)$ as a normalized intertwiner of the
tensor product of evaluation modules is equivalent to the following relations:
\vvn.3>
\begin{gather}
\label{relations-3}
R_{V_1,V_2}(x)\,v_1\ox v_2\,=\,v_1\ox v_2\,,\!
\\[5pt]
\bigl[\)R_{V_1,V_2}(x)\>, e_{a,b}\ox\id+\id\ox e_{a,b}\)\bigr]\,=\,0\,,\!
\notag
\\[4pt]
R_{V_1,V_2}(x)\>\bigl( x\,\id\ox e_{a,b}\>-
\tsum_{c=1}^N\>e_{a,c}\ox e_{c,b}\bigr)\,=\,
\bigl(x\,\id\ox e_{a,b}\>-\tsum_{c=1}^N e_{c,b}\ox e_{a,c}\bigr)\>
R_{V_1,V_2}(x)\,,\!
\notag
\end{gather}
for any $a\),b=1,\lc N$.
\vsk.2>
Let $V_1\),\)V_2\),\)V_3$ be Verma modules over $\gln$.
The corresponding \R-matrices satisfy the Yang-Baxter equation:
\vvn.3>
$$
R_{V_1,V_2}^{(1,2)}(x-y)\>R_{V_1,V_3}^{(1,3)}(x)\>R_{V_2,V_3}^{(2,3)}(y)\,=\,
R_{V_2,V_3}^{(2,3)}(y)\>R_{V_1,V_3}^{(1,3)}(x)\>R_{V_1,V_2}^{(1,2)}(x-y)\,.
\vv.2>
$$
\par
The formulated facts on the Yangian are well known;
for example, see~\cite{MNO}\).

\goodbreak
\subsection{The \qKZ/ and dynamical operators associated with $\gln\}$}
Let $V_1\lc V_n$ be Verma modules over $\gln$. Let $R_{V_i,V_j}(x)$ be
the corresponding rational \R-matrices. Let $p$, $\lan$ be nonzero complex
numbers. Denote by $T_u$ the difference operator acting on a function $f(u)$
by the formula
\vvn-.1>
$$
(T_uf)(u)\,=\,f(u+p)\,.
$$
\vsk.4>
Define the operators $K_1, \dots ,K_n$ acting on $\Vox$:
\vvn.2>
\begin{align*}
K_m(z\);\la)\,=\,\bigl(R_{V_1,V_m}^{(1,m)}(z_1\]-z_m\]-p)\ldots
R_{V_{m-1},V_m}^{(m-1,m)}(z_{m-1}\]-z_m\]-p)\bigr)\vpb{-1}\x{}\!\} &
\\[4pt]
{}\x\,\prod^N_{a=1}\,\la_a^{e_{a,a}^{(m)}}\,R^{(m,n)}_{V_m,V_n}(z_m\]-z_n)
\ldots R^{(m,m+1)}_{V_m,V_{m+1}}(z_m\]-z_{m+1})\,.&
\end{align*}
Introduce the difference operators $Z_1\lc Z_n$,
\ ${Z_i=\>\bigl(K_i(z;\la)\bigr)\vpb{-1}\>T_{z_i}}$, called
the \qKZ/ operators. They act on \$\Vox\)$-valued functions of $\zn\),\,\lan$.

\begin{thm}[\)\cite{FR}\)]
The \qKZ/ operators $Z_1\lc Z_n$ pairwise commute. In other words,
\vvn.2>
\begin{align*}
K_l(z_1\lc{}& z_m\]+p\lc z_n;\la)\>K_m(\zn;\la)\,={}
\\[4pt]
{}=\,K_m( &z_1\lc z_l\]+p\lc z_n; \la)\>K_l(\zn; \la)
\end{align*}
for all $\)l,m =1,\ldots ,n$.
\end{thm}

For $a\),b=1\lc N$, the element $e_{a,b}$ acts on $\Vox$ as
$\sum_{i=1}^n e_{a,b}^{(i)}$. Define the operators $L_1\lc L_N$
acting on $\Vox$:
\vvn.1>
$$
L_a(z\);\la)\,=\;\frac{(e_{a,a})^2}2\,-\,\sum_{i=1}^n\,z_i\)e_{a,a}^{(i)}\>-\)
\sum_{b=1}^N \sum_{1\le i<j\le n}\!\}e_{a,b}^{(i)}\>e_{b,a}^{(j)}\,-\>
\sum_{\fratop{b=1}{\smash{b\ne a}}}^N\,\frac{\la_b}{\la_a\]-\la_b}\,
(e_{a,b}\>e_{b,a}-e_{a,a})\,.
\vv-.2>
$$
Introduce the differential operators $D_1\lc D_N$,
\vvn.2>
$$
D_a\>=\,p\>\la_a\frac\der{\der\la_a}\,+\,L_a(z\);\la)\,,
\vv.1>
$$
called the dynamical operators. They act on \$\Vox\)$-valued
functions of $\zn\),\,\lan$.

\begin{thm}[\)\cite{TV5}\)]
\label{flat}
The dynamical operators $D_1\lc D_N$ pairwise commute.
\end{thm}

The theorem is proved by direct verification.

\subsection {Compatibility of \qKZ/ and dynamical operators}

\begin{thm}\label{compatibility}
The \qKZ/ operators $Z_1\lc Z_n$ commute with the dynamical operators
$D_1\lc D_N$.
\end{thm}

\vsk-.5>\vsk0>
\begin{proof}
Since the operators $Z_1\lc Z_n$ are invertible and pairwise commute,
the claim of the theorem is equivalent to the following identities:
\vvn.3>
$$
[\)Z_1\ldots Z_i\>,D_a\)]\,=0
\vv.2>
$$
for all $a=1\lc N$ and $i=1\lc n$. To simplify notations
we consider the case $i=1$. The general case is similar.
\goodbreak
\vsk.1>
A simple transformation converts identity $[\)Z_1\>,D_a\)]\)=\)0$ into identity
\vvn.2>
$$
\bigl[\)K_1(z\);\la)\>,L_a(z\);\la)\)]\,=\,0\,.
\vv-.4>
$$
Here \;$\dsize K_1(z\);\la)\,=\,\prod_{a=1}^N\,\la_a^{e_{a,a}^{(1)}}\,
R^{(1,n)}_{V_1,V_n}(z_1\]-z_n)\ldots R^{(1,2)}_{V_1,V_2}(z_1\]-z_2)$.
\,We have
\vvn-.1>
\begin{align}
L_a(z\);\la)\,=\;\frac{(e_{a,a})^2}2\,-\,
\sum_{\fratop{b=1}{\smash{b\ne a}}}^N\,\frac{\la_b}{\la_a\]-\la_b}\,
(e_{a,b}\>e_{b,a}-e_{a,a})\>-\>z_1\)e_{a,a}\)+{}
\label{L_a}
\\
{}+\,\sum_{i=2}^n\,(z_1\]-z_i)\)e_{a,a}^{(i)}\,-\,
\sum_{b=1}^N \sum_{1\le i<j\le n}\!\}e_{a,b}^{(i)}\>e_{b,a}^{(j)}\,&.
\notag
\end{align}
The term $z_1\)e_{a,a}$ commutes with $K_1(z\);\la)$. The first two terms
in the right hand side of \Ref{L_a} commute with the product
$R^{(1,n)}_{V_1,V_n}(z_1\]-z_n)\ldots R^{(1,2)}_{V_1,V_2}(z_1\]-z_2)$ since the
\R-matrices commute with coproducts. To proceed with the last two terms we are
using the second and third relations in \Ref{relations-3}\), and commutativity
of $R^{(1,i)}_{V_1,V_i}$ and $e_{c,d}^{(j)}$ for distinct $i,j$ not equal to
$1$. For example,
\vvn-.4>
\begin{alignat*}2
R^{(1,2)}_{V_1,V_2}(z_1\]-z_2) & \Bigl((z_1\]-z_2)\)e_{a,a}^{(2)}\,-\>
\sum_{b=1}^N\,e_{a,b}^{(1)}\>e_{b,a}^{(2)}\>+\>
\sum_{i=3}^n\,(z_1\]-z_i)\)e_{a,a}^{(i)}\>-{} &&
\\
&& \Llap{{}-\>
\sum_{b=1}^N\,\sum_{j=3}^n\,(e_{a,b}^{(1)}+e_{a,b}^{(2)})\>e_{b,a}^{(j)}\>-\>
\sum_{b=1}^N \sum_{3\le i<j\le n}\!\}e_{a,b}^{(i)}\>e_{b,a}^{(j)}\Bigr)}&
{}\>={}
\\[4pt]
{}={}\>& \Bigl((z_1\]-z_2)\)e_{a,a}^{(2)}\,-\>
\sum_{b=1}^N\,e_{b,a}^{(1)}\>e_{a,b}^{(2)}\>+\>
\sum_{i=3}^n\,(z_1\]-z_i)\)e_{a,a}^{(i)}\>-{}
\\
&& \Llap{{}-\>
\sum_{b=1}^N\,\sum_{j=3}^n\,(e_{a,b}^{(1)}+e_{a,b}^{(2)})\>e_{b,a}^{(j)}\>-\>
\sum_{b=1}^N \sum_{3\le i<j\le n}\!\}e_{a,b}^{(i)}\>e_{b,a}^{(j)}\Bigr)}&
\>R^{(1,2)}_{V_1,V_2}(z_1\]-z_2)\,,
\end{alignat*}
and consecutive calculations with other \R-matrices $R^{(1,i)}_{V_1,V_i}$ are
\vvn.3>
similar. Finally, we obtain
\begin{alignat}2
\label{last}
K_1(z\);\la)\>L_a(z\);\la)\, &{}=\,-\)z_1\)e_{a,a}\>K_1(z\);\la)\,+{}
\\[5pt]
&{}+\,\prod_{c=1}^N\,\la_c^{e_{c,c}^{(1)}}
\biggl(\frac{(e_{a,a})^2}2\,-\,
\sum_{\fratop{b=1}{\smash{b\ne a}}}^N\,\frac{\la_b}{\la_a\]-\la_b}\,
(e_{a,b}\>e_{b,a}-e_{a,a})\>+{}
\notag
\\
& {}\)+\,\sum_{i=2}^n\,(z_1\]-z_i)\)e_{a,a}^{(i)}\,-\,
\sum_{b=1}^N\,\sum_{j=2}^n\,e_{b,a}^{(1)}\>e_{a,b}^{(j)}\,-\,
\sum_{b=1}^N \sum_{2\le i<j\le n}\!\}e_{a,b}^{(i)}\>e_{b,a}^{(j)}\biggr)\x{}
\notag
\\[6pt]
&& \Llap{{}\x\,R^{(1,n)}_{V_1,V_n}(z_1\]-z_n)\ldots
R^{(1,2)}_{V_1,V_2}(z_1\]-z_2)\,.\!}
\notag
\end{alignat}
The right hand side of \Ref{last} can be transformed into
$L_a(z\);\la)\>K_1(z\);\la)$ in a straightforward way.
Theorem \ref{compatibility} is proved.
\end{proof}

\goodbreak
\subsection{Dynamical operators and the Weyl group}
The symmetric group $S_N$ is the Weyl group for the Lie algebra $\gln$.
There is a natural action, denoted by $\pi$, of $S_N$ on $U(\gln)$:
$\pi(w)(e_{a,b})=e_{w(a),w(b)}$ for $w\in S_N$ and $a\),b=1\lc N$.
\vsk.2>
Consider the space of linear differential operators of the form
\vvn.4>
$$
D\,=\,p\>\la_a\frac\der{\der\la_a}\,+\,M(z\);\la)\,,
\vv.4>
$$
where $a\in\lb 1\lc N\)\rb$ and $M$ is a \$U(\gln)$-valued function of
$\zn\),\,\lan$. The group $S_N$ acts on this space by the formula
\vvn.5>
$$
w\)D\,=\,p\>\la_{w(a)}\frac\der{\der\la_{w(a)}}\,+\,
\pi(w)\)M(z\);\la_{w(1)}\lc\la_{w(N)})\,.
\vv.3>
$$

\begin{thm}
We have \;$w\)D_a\)=\,D_{w(a)}$ \;for all $w\in S_N$ and $a=1\lc N$.
\end{thm}

The proof is straightfirward.

\end{document}